\documentclass[12pt,draft,leqno]{article}
\usepackage{amssymb, eucal, latexsym}

0in \oddsidemargin 1,2cm \evensidemargin 0in

\newtheorem{theorem}{Theorem}[section]
\newtheorem{lm}[theorem]{Lemma}
\newtheorem{exa}[theorem]{Example}

\newtheorem{pro}[theorem]{Proposition}
\newtheorem{defi}[theorem]{Definition}

\newtheorem{nota}[theorem]{Notation}
\newtheorem{notas}[theorem]{Notations}
\newtheorem{rem}[theorem]{Remark}

\newtheorem{fact}[theorem]{Fact}

\newtheorem{nist}[theorem]{}

\def\p{\varphi}
\def\a{\alpha}

\def\d{\delta}

\def\g{\gamma}
\def\GA{\Gamma}

\def\l{\lambda}

\def\s{\sigma}

\def\lra{\longrightarrow}

\def\sbe{\subseteq}
\def\spe{\supseteq}
\def\stm{\setminus}
\def\ems{\emptyset}
\def\nes{\neq\emptyset}

\def\cuk{\,\check{}\,}
\def\gek{\,\breve{}\,}

\def\unl{\underline}

\def\fa{\forall}

\def\we{\wedge}

\def\ap{^\prime}
\def\inv{^{-1}}
\def\st{\ |\ }

\def\ix{\infty_X}
\def\iy{\infty_Y}

\def\llx{\ll_{\rho}}
\def\lle{\ll_{\eta}}

\def\llcr{\ll_{C_\rho}}
\def\llce{\ll_{C_\eta}}

\def\nin{\not\in}

\def\card #1{\vert #1 \vert}

\def\CC{{\cal C}}

\def\TT{{\cal T}}

\def\Bo{{\bf Bool}}
\def\HC{{\bf HC}}

\def\VAL{{\bf VAL}}
\def\NAL{{\bf NAL}}
\def\PAL{{\bf PAL}}
\def\PALC{{\bf PALC}}
\def\PLC{{\bf PLC}}
\def\PLCC{{\bf PLCC}}

\def\VAL{{\bf DVAL}}

\def\2{\mbox{{\bf 2}}}
\def\3{\mbox{{\bf 3}}}

\def\int{\mbox{{\rm int}}}

\def\cl{\mbox{{\rm cl}}}
\def\CL{\mbox{{\rm Clust}}}
\def\BClu{\mbox{{\rm BClust}}}

\def\doc{\hspace{-1cm}{\em Proof.}~~}
\def\sq{\hspace*{\fill} \hbox{\vrule\vbox{\hrule\phantom{o}\hrule}\vrule}}
\def\sqs{\sq \vspace{2mm}}

\def\BBBB{{\rm I}\!{\rm B}}

\def\eset{\emptyset}   
\def\neset{\neq\emptyset} 

\title{{\LARGE\bf
A Generalization of De Vries Duality Theorem}\thanks{This paper
was supported by the project no. 26/2006 $``$Categorical Topology"
of
the Sofia University $``$St. Kl. Ohridski".}\\
\vspace{0.35cm}
{\large\bf Georgi Dimov}}

\author{}

\date{}

\begin{document}
\maketitle
\begin{abstract}
{\footnotesize
\noindent Generalizing Duality Theorem of H. de Vries, we define a
category which is dually equivalent to the category of all locally
compact Hausdorff spaces and all perfect maps between them.}
\end{abstract}

{\footnotesize {\em 2000 MSC:} primary 18A40, 54D45; secondary
06E15, 54C10, 54E05, 06E10.

{\em Keywords:}  Local contact algebra; Locally compact spaces;
Perfect maps; Duality.}

\footnotetext[1]{{\footnotesize {\em E-mail address:}
gdimov@fmi.uni-sofia.bg}}

\baselineskip = \normalbaselineskip

\section*{Introduction}

According to the famous Stone Duality Theorem (\cite{ST}), the
category of all zero-dimensi\-o\-nal compact Hausdorff spaces and
all continuous maps between them is dually equivalent to the
category $\Bo$ of all Boolean algebras and all Boolean
homomorphisms between them. In 1962, H. de Vries \cite{dV}
introduced the notion of {\em compingent Boolean algebra}\/ and
proved that the category of all compact Hausdorff spaces and all
continuous maps between them is dually equivalent to the category
of all complete compingent Boolean algebras and appropriate
morphisms between them. In 1997, Roeper \cite{R} defined the
notion of {\em region-based topology}\/  as one of the possible
formalizations of the ideas of De Laguna \cite{dL} and Whitehead
\cite{W} for a region-based theory of space.  Following
\cite{VDDB, DV}, the region-based topologies of Roeper appear here
as {\em local contact algebras}\/ (briefly, LCAs), because the
axioms which they satisfy almost coincide with the axioms of local
proxi\-mities of Leader \cite{LE}. In his paper \cite{R}, Roeper
proved the following theorem: there is a bijective correspondence
between all (up to homeomorphism) locally compact Hausdorff spaces
and all (up to isomorphism) complete LCAs. It generalizes the
theorem of de Vries \cite{dV} that there exists a bijective
correspondence between all (up to homeomorphism) compact Hausdorff
spaces and all (up to isomorphism) complete compingent Boolean
algebras. Here, using  Roeper's Theorem and the results of de
Vries \cite{dV}, a category dually equivalent to the category of
all locally compact Hausdorff spaces and all perfect maps between
them is defined (see Theorem \ref{gendv} bellow), genera\-lizing
in this way the  Duality Theorem of H. de Vries.

Let us mention that, using de Vries Duality Theorem, V. V.
Fedorchuk \cite{F} showed that the category  of all compact
Hausdorff spaces and all quasi-open maps between them is dually
equivalent to the category  of all complete compingent Bool\-ean
algebras and all complete Boolean homomorphisms between them
satisfying one simple condition, and that in \cite{D1, D2} some
extensions of the Fedorchuk Duality Theorem (\cite{F}) to some
categories whose objects are all locally compact Hausdorff spaces
are obtained.

We now fix the notations.

If $\CC$ denotes a category, we write $X\in \card\CC$ if $X$ is
 an object of $\CC$, and $f\in \CC(X,Y)$ if $f$ is a morphism of
 $\CC$ with domain $X$ and codomain $Y$.

All lattices are with top (= unit) and bottom (= zero) elements,
denoted respectively by 1 and 0. We do not require the elements
$0$ and $1$ to be distinct.

If $(X,\tau)$ is a topological space and $M$ is a subset of $X$,
we denote by $\cl_{(X,\tau)}(M)$ (or simply by $\cl(M)$ or
$\cl_X(M)$) the closure of $M$ in $(X,\tau)$ and by
$\int_{(X,\tau)}(M)$ (or briefly by $\int(M)$ or $\int_X(M)$) the
interior of $M$ in $(X,\tau)$. The Alexandroff compactification of
a locally compact Hausdorff non-compact space $X$ will be denoted
by $\a X$ and the added point by $\ix$ (i.e. $\a X=X\cup\{\ix\}$).

The  closed maps  between topological spaces are assumed to be
continuous but are not assumed to be onto. Recall that a map is
{\em perfect}\/ if it is closed and compact (i.e. point inverses
are compact sets).

\section{Preliminaries}
%

\begin{defi}\label{conalg}
\rm
An algebraic system $\underline {B}=(B,0,1,\vee,\we, {}^*, C)$ is
called a {\it contact algebra} (abbreviated as CA)
 if
$(B,0,1,\vee,\we, {}^*)$ is a Boolean algebra (where the operation
$``$complement" is denoted by $``\ {}^*\ $")
  and $C$
is a binary relation on $B$, satisfying the following axioms:

\smallskip

\noindent (C1) If $a\not= 0$ then $aCa$;\\
(C2) If $aCb$ then $a\not=0$ and $b\not=0$;\\
(C3) $aCb$ implies $bCa$;\\
(C4) $aC(b\vee c)$ iff $aCb$ or $aCc$.

\smallskip

\noindent Usually, we shall simply write $(B,C)$ for a contact
algebra. The relation $C$  is called a {\em  contact relation}.
When $B$ is a complete Boolean algebra, we will say that $(B,C)$
is a {\em complete contact algebra}\/ (abbreviated as CCA).

We will say that two CA's $(B_1,C_1)$ and $(B_2,C_2)$ are  {\em
CA-isomorphic} iff there exists a Boolean isomorphism $\p:B_1\lra
B_2$ such that, for each $a,b\in B_1$, $aC_1 b$ iff $\p(a)C_2
\p(b)$. Note that in this paper, by a $``$Boolean isomorphism" we
understand an isomorphism in the category $\Bo$.

A CA $(B,C)$  is called {\em connected}\/ if it satisfies the
following axiom:

\smallskip

\noindent (CON) If $a\neq 0,1$ then $aCa^*$.

\smallskip

A contact algebra $(B,C)$ is called a {\it  normal contact
algebra} (abbreviated as NCA) (\cite{dV,F}) if it satisfies the
following axioms (we will write $``-C$" for $``not\ C$"):

\smallskip

\noindent (C5) If $a(-C)b$ then $a(-C)c$ and $b(-C)c^*$ for some $c\in B$;\\
(C6) If $a\not= 1$ then there exists $b\not= 0$ such that
$b(-C)a$.

\smallskip

\noindent A normal CA is called a {\em complete normal contact
algebra} (abbreviated as CNCA) if it is a CCA. The notion of
normal contact algebra was introduced by Fedorchuk \cite{F} under
the name {\em Boolean $\d$-algebra}\/ as an equivalent expression
of the notion of compingent Boolean algebra of de Vries. We call
such algebras $``$normal contact algebras" because they form a
subclass of the class of contact algebras.

Note that if $0\neq 1$ then the axiom (C2) follows from the axioms
(C6) and (C4).

For any CA $(B,C)$, we define a binary relation  $``\ll_C $"  on
$B$ (called {\em non-tangential inclusion})  by $``\ a \ll_C b
\leftrightarrow a(-C)b^*\ $". Sometimes we will write simply
$``\ll$" instead of $``\ll_C$".
\end{defi}

The relations $C$ and $\ll$ are inter-definable. For example,
normal contact algebras could be equivalently defined (and exactly
in this way they were defined (under the name of compingent
Boolean algebras) by de Vries in \cite{dV}) as a pair of a Boolean
algebra $B=(B,0,1,\vee,\we,{}^*)$ and a binary relation $\ll$ on
$B$ subject to the following axioms:

\smallskip

\noindent ($\ll$1) $a\ll b$ implies $a\leq b$;\\
($\ll$2) $0\ll 0$;\\
($\ll$3) $a\leq b\ll c\leq t$ implies $a\ll t$;\\
($\ll$4) $a\ll c$ and $b\ll c$ implies $a\vee b\ll c$;\\
($\ll$5) If  $a\ll c$ then $a\ll b\ll c$  for some $b\in B$;\\
($\ll$6) If $a\neq 0$ then there exists $b\neq 0$ such that $b\ll
a$;\\
($\ll$7) $a\ll b$ implies $b^*\ll a^*$.

\smallskip

Note that if $0\neq 1$ then the axiom ($\ll$2) follows from the
axioms ($\ll$3), ($\ll$4), ($\ll$6) and ($\ll$7).

\smallskip

Obviously, contact algebras could be equivalently defined as a
pair of a Boolean algebra $B$ and a binary relation $\ll$ on $B$
subject to the  axioms ($\ll$1)-($\ll$4) and ($\ll$7).

\smallskip

It is easy to see that axiom (C5) (resp., (C6)) can be stated
equivalently in the form of ($\ll$5) (resp., ($\ll$6)).

\begin{exa}\label{extrcr}
\rm Let $B$ be a Boolean algebra. Then there exist the largest and
the smallest contact relations on $B$; the largest one, $\rho_l$,
is defined by $a\rho_l b$ iff $a\neq 0$ and $b\neq 0$, and the
smallest one, $\rho_s$, by $a\rho_s b$ iff $a\wedge b\neq 0$.

Note that, for $a,b\in B$, $a\ll_{\rho_s} b$ iff $a\le b$; hence
$a\ll_{\rho_s} a$, for any $a\in B$. Thus $(B,\rho_s)$ is a normal
contact algebra.
\end{exa}

\begin{exa}\label{rct}
\rm Recall that a subset $F$ of a topological space $(X,\tau)$ is
called {\em regular closed}\/ if $F=\cl(\int (F))$. Clearly, $F$
is regular closed iff it is the closure of an open set.

For any topological space $(X,\tau)$, the collection $RC(X,\tau)$
(we will often write simply $RC(X)$) of all regular closed subsets
of $(X,\tau)$ becomes a complete Boolean algebra
$(RC(X,\tau),0,1,\we,\vee,{}^*)$ under the following operations:
$$ 1 = X,  0 = \emptyset, F^* = \cl(X\stm F), F\vee G=F\cup G,
F\we G =\cl(\int(F\cap G)).
$$
The infinite operations are given by the following formulas:
$\bigvee\{F_\g\st \g\in\GA\}=\cl(\bigcup\{F_\g\st
\g\in\GA\})(=\cl(\bigcup\{\int(F_\g)\st \g\in\GA\})),$ and
$\bigwedge\{F_\g\st \g\in\GA\}=\cl(\int(\bigcap\{F_\g\st
\g\in\GA\})).$

It is easy to see that setting $F \rho_{(X,\tau)} G$ iff $F\cap
G\not = \ems$, we define a contact relation $\rho_{(X,\tau)}$ on
$RC(X,\tau)$; it is called a {\em standard contact relation}. So,
$(RC(X,\tau),\rho_{(X,\tau)})$ is a CCA (it is called a {\em
standard contact algebra}). We will often write simply $\rho_X$
instead of $\rho_{(X,\tau)}$. Note that, for $F,G\in RC(X)$,
$F\ll_{\rho_X}G$ iff $F\sbe\int_X(G)$.

Clearly, if $(X,\tau)$ is a normal Hausdorff space then the
standard contact algebra $(RC(X,\tau),\rho_{(X,\tau)})$ is a
complete NCA.

A subset $U$ of $(X,\tau)$ such that $U=\int(\cl(U))$ is said to
be {\em regular open}. The set of all regular open subsets of
$(X,\tau)$ will be denoted by $RO(X,\tau)$ (or briefly, by
$RO(X)$). Define Boolean operations and contact $\d_X$ in $RO(X)$
as follows: $U\vee V=\int(\cl(U\cup V))$, $U\wedge V=U\cap V$,
$U^{*}=\int(X\setminus U)$, $0=\eset$, $1=X$ and $U\d_X V$ iff
$\cl (U)\cap \cl (V)\neset$. Then $(RO(X),\d_X)$ is a CA. This
algebra is also complete, considering the infinite meet
$\bigwedge\{U_i\st i\in I\} =\int(\bigcap_{i\in I}U_{i})$.

Note that $(RO(X),\d_X)$ and $(RC(X),\rho_X)$ are isomorphic CAs.
The isomorphism $f$ between them is defined by  $f(U)=\cl(U)$, for
every $U\in RO(X)$.
\end{exa}

The following notion is a lattice-theoretical counterpart of the
corresponding notion from the theory of proximity spaces (see
\cite{NW}):

\begin{nist}\label{defcluclan}
\rm Let $(B,C)$ be a CA. Then  a non-empty subset $\s $ of $B$ is
called a {\em cluster in} $(B,C)$
if the
following conditions are satisfied:

\smallskip

\noindent (K1) If $a,b\in\s $ then $aCb$;\\
(K2) If $a\vee b\in\s $ then $a\in\s $ or $b\in\s $;\\
(K3) If $aCb$ for every $b\in\s $, then $a\in\s $.

\smallskip

\noindent The set of all clusters in $(B,C)$ will be denoted
denoted by $\CL(B,C)$.
\end{nist}

The next  assertion can be proved exactly as Lemma 5.6 of
\cite{NW}:

\begin{fact}\label{fact29}
If $\s_1,\s_2$ are two clusters in a normal contact algebra
$(B,C)$ and $\s_1\sbe \s_2$ then $\s_1=\s_2$.
\end{fact}

\begin{fact}\label{confact}{\rm (\cite{BG})}
Let $(X,\tau)$ be a topological space. Then the standard contact
algebra $(RC(X,\tau),\rho_{(X,\tau)})$ is connected iff the space
$(X,\tau)$ is connected.
\end{fact}

The following notion is a lattice-theoretical counterpart of the
Leader's notion of {\em local proximity} (\cite{LE}):

\begin{defi}\label{locono}{\rm (\cite{R})}
\rm An algebraic system $\underline {B}_{\, l}=(B,0,1,\vee,\we,
{}^*, \rho, \BBBB)$ is called a {\it local contact algebra}
(abbreviated as LCA)   if $(B,0,1, \vee,\we, {}^*)$ is a Boolean
algebra, $\rho$ is a binary relation on $B$ such that $(B,\rho)$
is a CA, and $\BBBB$ is an ideal (possibly non proper) of $B$,
satisfying the following axioms:

\smallskip

\noindent(BC1) If $a\in\BBBB$, $c\in B$ and $a\ll_\rho c$ then
$a\ll_\rho b\ll_\rho c$ for some $b\in\BBBB$  (see \ref{conalg}
for
$``\ll_\rho$");\\
(BC2) If $a\rho b$ then there exists an element $c$ of $\BBBB$
such that
$a\rho (c\we b)$;\\
(BC3) If $a\neq 0$ then there exists  $b\in\BBBB\stm\{0\}$ such
that $b\ll_\rho a$.

\smallskip

Usually, we shall simply write  $(B, \rho,\BBBB)$ for a local
contact algebra.  We will say that the elements of $\BBBB$ are
{\em bounded} and the elements of $B\stm \BBBB$  are  {\em
unbounded}. When $B$ is a complete Boolean algebra,  the LCA
$(B,\rho,\BBBB)$ is called a {\em complete local contact algebra}
(abbreviated as CLCA).

We will say that two local contact algebras $(B,\rho,\BBBB)$ and
$(B_1,\rho_1,\BBBB_1)$ are  {\em LCA-isomorphic} iff there exists
a Boolean isomorphism $\p:B\lra B_1$ such that, for $a,b\in B$,
$a\rho b$ iff $\p(a)\rho_1 \p(b)$, and $\p(a)\in\BBBB_1$ iff
$a\in\BBBB$.

An LCA $(B,\rho,\BBBB)$ is called {\em connected}\/ if the CA
$(B,\rho)$ is connected.
\end{defi}

\begin{rem}\label{conaln}
\rm Note that if $(B,\rho,\BBBB)$ is a local contact algebra and
$1\in\BBBB$ then $(B,\rho)$ is a normal contact algebra.
Conversely, any normal contact algebra $(B,C)$ can be regarded as
a local contact algebra of the form $(B,C,B)$.
\end{rem}

The following lemmas
are lattice-theoretical
counterparts of some theorems from Leader's paper \cite{LE}.

\begin{lm}\label{Alexprn}{\rm (\cite{VDDB})}
Let $(B,\rho,\BBBB)$ be a local contact algebra. Define a binary
relation $``C_\rho$" on $B$ by
\begin{equation}\label{crho}
aC_\rho b\ \mbox{ iff }\ a\rho b\ \mbox{ or }\ a,b\not\in\BBBB.
\end{equation}
Then $``C_\rho$", called the\/ {\em Alexandroff extension of}\/
$\rho$, is a normal contact relation on $B$ and $(B,C_\rho)$ is a
normal contact algebra.
\end{lm}

\begin{lm}\label{neogrn}{\rm (\cite{VDDB})}
Let $\underline {B}_{\, l}=(B,\rho,\BBBB)$ be a local contact
algebra and let $1\not\in\BBBB$. Then $\s_\infty^{\underline
{B}_{\, l}}=\{b\in B\st b\not\in\BBBB\}$ is a cluster in
$(B,C_\rho)$ (see \ref{Alexprn} for the notation $``C_\rho$").
(Sometimes we will simply write  $\s_\infty$
instead of $\ \s_\infty^{\underline {B}_{\, l}}$.)
\end{lm}

\begin{defi}\label{boundcl}
\rm Let $(B,\rho,\BBBB)$ be a local contact algebra. A cluster
$\s$ in $(B,C_\rho)$ (see \ref{Alexprn}) is called {\em bounded}\/
if $\s\cap\BBBB\nes$. The set of all bounded clusters in
$(B,C_\rho)$ will be denoted by $\BClu(B,\rho,\BBBB)$.
\end{defi}

\begin{fact}\label{bstar}
Let $(B,\rho,\BBBB)$ be a local contact algebra and  $\s$ be a
bounded cluster in $(B,C_\rho)$ (see \ref{Alexprn}). Then there
exists  $b\in \BBBB$  such that  $b^*\nin\s$.
\end{fact}

\doc
Let $b_0\in\s\cap\BBBB$. Since $b_0\ll_\rho 1$, (BC1) implies that
there exists $b\in\BBBB$ such that $b_0\ll_\rho b$. Then
$b_0(-\rho)b^*$ and since $b_0\in\BBBB$, we obtain that
$b_0(-C_\rho)b^*$. Thus $b^*\nin\s$. \sqs

\begin{nota}\label{compregn}
\rm Let $(X,\tau)$ be a topological space. We denote by
$CR(X,\tau)$ the family of all compact regular closed subsets of
$(X,\tau)$. We will often write  $CR(X)$ instead of $CR(X,\tau)$.

 If $x\in X$ then we
set:
\begin{equation}\label{sxvx}
\s_x=\{F\in RC(X)\st x\in F\}.
\end{equation}
\end{nota}

\begin{fact}\label{stanlocn}
Let $(X,\tau)$ be a locally compact Hausdorff space. Then the
triple
$$(RC(X,\tau),\rho_{(X,\tau)}, CR(X,\tau))$$
 (see \ref{rct} for $\rho_{(X,\tau)}$)
  is a complete local contact algebra   (\cite{R}). It is called a
{\em standard local contact algebra}.

For every $x\in X$, $\s_x$ is a bounded cluster in
$(RC(X),C_{\rho_X})$ (see (\ref{sxvx}) and (\ref{crho}) for the
notations).
\end{fact}

We will need a lemma from \cite{CNG}:

\begin{lm}\label{isombool}
Let $X$ be a dense subspace of a topological space $Y$. Then the
functions $r_{X,Y}:RC(Y)\lra RC(X)$, $F\lra F\cap X$, and
$e_{X,Y}:RC(X)\lra RC(Y)$, $G\lra \cl_Y(G)$, are Boolean
isomorphisms between Boolean algebras $RC(X)$ and $RC(Y)$, and
$e_{X,Y}\circ r_{X,Y}=id_{RC(Y)}$, $r_{X,Y}\circ
e_{X,Y}=id_{RC(X)}$. (We will often write $r_X, e_X$ instead of
$r_{X,Y},e_{X,Y}$, respectively.)
\end{lm}

The next proposition is well known (see, e.g., \cite{AP}):

\begin{pro}\label{perfect1}
Let $f:X\lra Y$ be a perfect map between two locally compact
Hausdorff non-compact spaces. Then the map $f$ has a continuous
extension $\a(f):\a X\lra\a Y$; moreover, $\a(f)(\ix)=\iy$.
\end{pro}

For all undefined here notions and notations see \cite{AHS, J, E,
NW, Si}.

\section{The Results}

The next  theorem was proved by Roeper \cite{R}.  We will give a
sketch of its proof; it follows the plan of the proof presented in
\cite{VDDB}. The notations and the facts stated here will be used
later on.

\begin{theorem}\label{roeperl}{\rm (P. Roeper \cite{R})}
There exists a bijective correspondence between the class of all
(up to isomorphism) complete local contact algebras and the class
of all (up to homeomorphism) locally compact Hausdorff spaces.
\end{theorem}

\noindent{\em Sketch of the Proof.}~ (A) Let $(X,\tau)$ be a
locally compact Hausdorff space. We put
\begin{equation}\label{psit1}
\Psi^t(X,\tau)=(RC(X,\tau),\rho_{(X,\tau)},CR(X,\tau))
\end{equation}
(see \ref{stanlocn} and \ref{compregn} for the notations).

\noindent(B)~ Let $\unl{B}_{\, l}=(B,\rho,\BBBB)$ be a complete
local contact algebra. Let $C=C_\rho$ be the Alexandroff extension
of $\rho$ (see \ref{Alexprn}). Then, by  \ref{Alexprn}, $(B,C)$ is
a complete normal contact algebra. Put $X=\CL(B,C)$ and let $\TT$
be the topology on $X$ having as a closed base the family
$\{\l_{(B,C)}(a)\st a\in B\}$ where, for every $a\in B$,
\begin{equation}\label{h}
\l_{(B,C)}(a) = \{\s \in X\st  a \in \s\}.
\end{equation}
Sometimes we will write simply $\l_B$ instead of $\l_{(B,C)}$.

It can be proved that $(X,\TT)$ is a compact Hausdorff space and
\begin{equation}\label{isom}
\l_B:(B,C)\lra (RC(X),\rho_X) \mbox{ is a CA-isomorphism.}
\end{equation}

\smallskip

\noindent(B1)~ Let $1\in\BBBB$. Then $C=\rho$ and $\BBBB=B$, so
that $(B,\rho,\BBBB)=(B,C,B)=(B,C)$ is a complete normal contact
algebra (see \ref{conaln}), and we put
\begin{equation}\label{phiapcn}
\Psi^a(B,\rho,\BBBB)=\Psi^a(B,C,B)=\Psi^a(B,C)=(X,\TT).
\end{equation}

\medskip

\noindent(B2)~ Let $1\not\in\BBBB$. Then, by Lemma \ref{neogrn},
the set $\s_\infty=\{b\in B\st b\not\in\BBBB\}$ is a cluster in
$(B,C)$ and, hence, $\s_\infty\in X$.  Let $L=X\stm\{\s_\infty\}$.
Then
\begin{equation}\label{L}
L=\BClu(B,\rho,\BBBB),
\end{equation}
i.e.  $L$ is the set of all bounded clusters of  $(B,C_\rho)$
(sometimes we will write $L_{\unl{B}_{\, l}}$ or $L_B$ instead of
$L$);
 let the topology $\tau(=\tau_{\unl{B}_{\, l}})$ on $L$ be the
subspace topology, i.e. $\tau=\TT|_L $. Then $(L,\tau)$ is a
locally compact Hausdorff space. We put
\begin{equation}\label {phiapc}
\Psi^a(B,\rho,\BBBB)=(L,\tau).
\end{equation}

Let
%
$\l^l_B(a)=\l_B(a)\cap L$,
%
for each $a\in B$. One can show that $X=\a L$ and
\begin{equation}\label{hapisom}
\l^l_B: (B,\rho,\BBBB)\lra (RC(L),\rho_L, CR(L)) \mbox{ is an
LCA-isomorphism.}
\end{equation}

\medskip

\noindent(C)~ For every CLCA $(B,\rho,\BBBB)$ and every $a\in B$,
set
\begin{equation}\label{lbg}
\l^g_B(a)=\l_B(a)\cap\Psi^a(B,\rho,\BBBB).
\end{equation}
Then, by (\ref{isom}) and (\ref{hapisom}), we get that
\begin{equation}\label{hapisomn}
\l^g_B: (B,\rho,\BBBB)\lra (\Psi^t\circ\Psi^a)(B,\rho,\BBBB)
\mbox{ is an LCA-isomorphism.}
\end{equation}

\medskip

\noindent(D)~~ Let  $(Y,\tau)$ be a locally compact Hausdorff
space.   It can be shown that the map
\begin{equation}\label{homeo}
t_{(Y,\tau)}:(Y,\tau)\lra\Psi^a(\Psi^t(Y,\tau)),
\end{equation}
defined by  $t_{(Y,\tau)}(y)=\{F\in RC(Y,\tau)\st y\in
F\}(=\s_y)$, for every $y\in Y$, is a homeomorphism;  we will
often write simply $t_Y$ instead of $t_{(Y,\tau)}$.

Therefore $\Psi^a(\Psi^t(Y,\tau))$ is homeomorphic to $(Y,\tau)$
and $\Psi^t(\Psi^a(B,\rho,\BBBB))$ is LCA-isomorphic to
$(B,\rho,\BBBB)$.
 \sqs

\begin{defi}\label{dval}{\rm (De Vries \cite{dV})}
\rm Let  $\HC$ be the category of all compact Hausdorff spaces and
all continuous maps between them.

Let $\VAL$ be the category whose objects are all complete NCAs and
whose morphisms are all functions $\p:(A,C)\lra (B,C\ap)$ between
the objects of $\VAL$ satisfying the
conditions:\\
(DVAL1) $\p(0)=0$;\\
(DVAL2) $\p(a\we b)=\p(a)\we \p(b)$, for all $a,b\in A$;\\
(DVAL3) If $a, b\in A$ and $a\ll_C b$ then $(\p(a^*))^*\ll_{C\ap}
\p(b)$;\\
(DVAL4) $\p(a)=\bigvee\{\p(b)\st b\ll_{C} a\}$, for every $a\in
A$,

\medskip

{\noindent}and let the composition $``\ast$" of two morphisms
$\p_1:(A_1,C_1)\lra (A_2,C_2)$ and $\p_2:(A_2,C_2)\lra (A_3,C_3)$
of $\VAL$ be defined by the formula
\begin{equation}\label{diamc}
\p_2\ast\p_1 = (\p_2\circ\p_1)\gek,
\end{equation}
 where, for every
function $\psi:(A,C)\lra (B,C\ap)$ between two objects of $\VAL$,
$\psi\gek:(A,C)\lra (B,C\ap)$ is defined as follows:
\begin{equation}\label{cukfc}
\psi\gek(a)=\bigvee\{\psi(b)\st b\ll_{C} a\},
\end{equation}
for every $a\in A$.
\end{defi}

De Vries \cite{dV} proved the following duality theorem:

\begin{theorem}\label{dvth}
The categories $\HC$ and $\VAL$ are dually equivalent. In more
details, let $\Phi^t:\HC\lra\VAL$ be the contravariant functor
defined by $\Phi^t(X,\tau)=(RC(X,\tau),\rho_X)$, for every
$X\in\card\HC$, and $\Phi^t(f)(G)=\cl(f\inv(\int(G)))$, for every
$f\in\HC(X,Y)$ and every $G\in RC(Y)$, and let
$\Phi^a:\VAL\lra\HC$ be the contravariant functor defined by
$\Phi^a(A,C)=\Psi^a(A,C)$, for every $(A,C)\in\card\VAL$, and
$\Phi^a(\p)(\s\ap)=\{a\in A\st$if $b\ll_C a^*$ then
$(\p(b))^*\in\s\ap\}$, for every $\p\in\VAL((A,C),(B,C\ap))$ and
for every $\s\ap\in\CL(B,C\ap)$; then  $\l:
Id_{\,\VAL}\lra\Phi^t\circ\Phi^a$, where $\l(A,C)=\l_{(A,C)}$ (see
(\ref{h}) and (\ref{isom}) for the notation $\l_{(A,C)}$), for
every $(A,C)\in\card\VAL$, and
$t:Id_{\,\HC}\lra\Phi^a\circ\Phi^t$, where $t(X)=t_X$ (see
(\ref{homeo}) for the notation $t_X$), for every $X\in\card\HC$,
are natural isomorphisms.
\end{theorem}

In \cite{dV}, de Vries uses the regular open sets instead of
regular closed sets, as we do, so that we present here the
translations of his definitions for the case of regular closed
sets.

\begin{defi}\label{pal}
\rm We will denote by $\PLC$ the category of all locally compact
Hausdorff spaces and all perfect maps between them.

Let $\PAL$ be the category whose objects are all complete LCAs and
whose morphisms are all functions $\p:(A,\rho,\BBBB)\lra
(B,\eta,\BBBB\ap)$ between the objects of $\PAL$ satisfying the
conditions:\\
(PAL1) $\p(0)=0$;\\
(PAL2) $\p(a\we b)=\p(a)\we \p(b)$, for all $a,b\in A$;\\
(PAL3) If $a\in\BBBB, b\in A$ and $a\llx b$ then $(\p(a^*))^*\lle
\p(b)$;\\
(PAL4) For every $b\in\BBBB\ap$ there exists $a\in\BBBB$ such that
$b\le\p(a)$;\\
(PAL5) If $a\in\BBBB$ then $\p(a)\in\BBBB\ap$;\\
(PAL6) $\p(a)=\bigvee\{\p(b)\st b\ll_{C_\rho} a\}$, for every
$a\in A$ (see (\ref{crho}) for $C_\rho$);

\medskip

{\noindent}let the composition $``\diamond$" of two morphisms
$\p_1:(A_1,\rho_1,\BBBB_1)\lra (A_2,\rho_2,\BBBB_2)$ and
$\p_2:(A_2,\rho_2,\BBBB_2)\lra (A_3,\rho_3,\BBBB_3)$ of $\PAL$ be
defined by the formula
\begin{equation}\label{diam}
\p_2\diamond\p_1 = (\p_2\circ\p_1)\cuk,
\end{equation}
 where, for every
function $\psi:(A,\rho,\BBBB)\lra (B,\eta,\BBBB\ap)$ between two
objects of $\PAL$, $\psi\cuk:(A,\rho,\BBBB)\lra (B,\eta,\BBBB\ap)$
is defined as follows:
\begin{equation}\label{cukf}
\psi\cuk(a)=\bigvee\{\psi(b)\st b\ll_{C_\rho} a\},
\end{equation}
for every $a\in A$.

 By $\NAL$ we denote the full
subcategory of $\PAL$ having as objects all CNCAs (i.e., those
CLCAs $(A,\rho,\BBBB)$ for which $\BBBB=A$).
\end{defi}

Note that the categories $\VAL$ and $\NAL$ are isomorphic (it can
be even said that they are identical) because the axiom (PAL5) is
tri\-vially fulfilled in the category $\VAL$ (indeed, all elements
of its objects are bounded), the axiom (PAL4) follows immediately
from the obvious fact that $\p(1)=1$ for every $\VAL$-morphism
$\p$, and the compositions are the same.

 We will generalize the  Duality Theorem of de Vries showing that
the categories $\PAL$ and $\PLC$ are dually equivalent.

We will first show that $\PAL$ is indeed a category.

\begin{lm}\label{pf1}
Let us regard two functions $\p:(A,\rho,\BBBB)\lra
(B,\eta,\BBBB\ap)$ and $\psi:(B,\eta,\BBBB\ap)\lra
(B_1,\eta_1,\BBBB_1\ap)$
between  CLCAs. Then:\\
(a) If $\p$ satisfies condition (PAL2) then $\p$ is an order preserving function;\\
(b) If $\p$ satisfies conditions (PAL1) and
(PAL2) then  $\p(a^*)\le (\p(a))^*$, for every $a\in A$;\\
(c) Let $\p$ satisfy conditions (PAL3) and (PAL5).  If $a,b\in A$
and $a\ll_{C_\rho} b$ then $(\p(a^*))^*\ll_{C_\eta}\p(b)$. Hence,
if $\p$ satisfies in addition conditions (PAL1) and (PAL2) then
$\p(a)\ll_{C_\eta} \p(b)$;\\
(d) If $\p$ satisfies conditions (PAL1) and (PAL3) then  $\p(1)=1$;\\
(e) If $\p$ satisfies condition (PAL2) then $\p\cuk$ satisfies
conditions (PAL2) and (PAL6) (see (\ref{cukf}) for $\p\cuk$);\\
(f) If $\p$ satisfies condition (PAL6) then $\p=\p\cuk$;\\
(g) If $\p$ satisfies condition (PAL2) then
$(\p\cuk)\cuk=\p\cuk$;\\
(h) If $\p$ and $\psi$ satisfy condition (PAL2) and $\p$ satisfies
in addition conditions (PAL1), (PAL3) and (PAL5) then
$(\psi\circ\p)\cuk=(\psi\cuk\circ\p\cuk)\cuk$.
\end{lm}

\doc The properties (a), (b), (d) and (f) are clearly fulfilled, and (g)
follows from (e) and (f).

\medskip

 \noindent(c) Let $a,b\in A$ and $a\ll_{C_\rho}
b$. Then $a\llx b$ and at least one of the elements $a$ and $b^*$
is bounded.

Let $a\in\BBBB$. Then (PAL3) implies that $(\p(a^*))^*\lle\p(b)$.
By (BC1), there exists $c\in\BBBB$ such that $a\llx c$. Hence,
using again (PAL3), we get that $(\p(a^*))^*\lle\p(c)$. Since
$\p(c)\in\BBBB\ap$ (according to (PAL5)), we obtain that
$(\p(a^*))^*\in\BBBB\ap$. Therefore, $(\p(a^*))^*\llce\p(b)$.

Let now $b^*\in\BBBB$. Since $b^*\llcr a^*$, we get, by the
previous case, that $(\p(b))^*\llce\p(a^*)$. Thus
$(\p(a^*))^*\llce\p(b)$.

\medskip

\noindent(e) By (a), for every $a\in A$, $\p\cuk(a)\le\p(a)$. Let
$a\in A$. If $c\llcr a$ then there exists $d_c\in A$ such that
$c\llcr d_c\llcr a$; hence $\p(c)\le\p\cuk(d_c)$. Now,
$\p\cuk(a)=\bigvee\{\p(c)\st c\llcr a\}\le\bigvee\{\p\cuk(d_c)\st
c\llcr a\}\le\bigvee\{\p\cuk(e)\st e\llcr a\}\le\bigvee\{\p(e)\st
e\llcr a\}=\p\cuk(a)$. Thus, $\p\cuk(a)=\bigvee\{\p\cuk(e)\st
e\llcr a\}$. So, $\p\cuk$ satisfies (PAL6). Further, let $a,b\in
A$. Then $\p\cuk(a)\we\p\cuk(b)=\bigvee\{\p(d)\we\p(e)\st
d\ll_{C_\rho} a, e\ll_{C_\rho} b\}=\bigvee\{\p(d\we e)\st
d\ll_{C_\rho} a, e\ll_{C_\rho} b\}=\bigvee\{\p(c)\st c\ll_{C_\rho}
a\we b\}=\p\cuk(a\we b)$. So, (PAL2) is fulfilled.

\medskip

\noindent(h) Since $\p\cuk(a)\le \p(a)$ for every $a\in A$, and
$\psi\cuk (b)\le\psi(b)$ for every $b\in B$, we get that
$\psi\cuk(\p\cuk(a))\le\psi(\p(a))$, for every $a\in A$. Hence,
using (\ref{cukf}), we obtain that
$(\psi\cuk\circ\p\cuk)\cuk(a)\le(\psi\circ\p)\cuk(a)$, for every
$a\in A$. Further, by (\ref{cukf}), for every $a\in A$,
$(\psi\circ\p)\cuk(a)=\bigvee\{\psi(\p(e))\st e\llcr a\}$ and
$(\psi\cuk\circ\p\cuk)\cuk(a)=\bigvee\{\psi\cuk(\p\cuk(b))\st
b\llcr a\}=\bigvee\{\bigvee\{\psi(c)\st c\llce\p\cuk(b)\}\st
b\llcr a\}$. Let $a\in A$ and $e\llcr a$. Then there exist $b,d\in
A$ such that $e\llcr d\llcr b\llcr a$. Set $c=\p(e)$. Then, by
(c), $c\llce\p(d)\le\p\cuk(b)$. Hence
$\psi(\p(e))=\psi(c)\le(\psi\cuk\circ\p\cuk)\cuk(a)$. We conclude
that $(\psi\circ\p)\cuk(a)\le(\psi\cuk\circ\p\cuk)\cuk(a)$.
Therefore the desired equality is proved. \sqs

\begin{pro}\label{compos}
Let $\p_i:(A_i,\rho_i,\BBBB_i)\lra
(A_{i+1},\rho_{i+1},\BBBB_{i+1})$, where $i=1,2$,
be two
functions between CLCAs and let $\p_1$ and $\p_2$ satisfy
conditions (PAL1)-(PAL5). Then the function $\p_2\circ\p_1$
satisfies conditions (PAL1)-(PAL5).
\end{pro}

\doc Let $a\in\BBBB_1$,  $b\in A_1$ and $a\ll_{\rho_1} b$. Then,
by (BC1), there exists $c\in\BBBB_1$ such that $a\ll_{\rho_1}
c\ll_{\rho_1} b$. From (PAL3) we get that
$(\p_1(a^*))^*\ll_{\rho_2} \p_1(c)$. Then, since
$\p_1(c)\in\BBBB_2$ (by (PAL5)), $(\p_1(a^*))^*\in\BBBB_2$. Now,
using twice (PAL3), we obtain that
$(\p_1(a^*))^*\ll_{\rho_2}\p_1(b)$ and
$(\p_2(\p_1(a^*)))^*\ll_{\rho_3}\p_2(\p_1(b))$. Hence, the
function $\p_2\circ\p_1$ satisfies condition (PAL3). The rest is
obvious. \sqs

\begin{pro}\label{compos1}
Let $\p:(A,\rho,\BBBB)\lra (B,\eta,\BBBB\ap)$ be a function
between CLCAs and let $\p$  satisfies conditions (PAL1)-(PAL5).
Then the function $\p\cuk$ (see (\ref{cukf})) satisfies conditions
(PAL1)-(PAL6) (i.e., it is a $\PAL$-morphism).
\end{pro}

\doc Obviously, for every $a\in A$, $\p\cuk(a)\le\p(a)$. Hence,
$\p\cuk(0)=0$, i.e. (PAL1) is fulfilled. For (PAL2) and (PAL6) see
\ref{pf1}(e). Let $a\in\BBBB, b\in A$ and $a\llx b$. Then, by
(BC1), there exist $c,d\in\BBBB$ such that $a\llx c\llx d\llx b$.
Thus $a\llcr c\llcr d\llcr b$ and hence $c^*\llcr a^*$. We obtain
that $\p(d)\le\p\cuk(b)$ and $\p(c^*)\le\p\cuk(a^*)$. Hence
$(\p\cuk(a^*))^*\le (\p(c^*))^*\lle \p(d)\le\p\cuk(b)$. Therefore,
$(\p\cuk(a^*))^*\lle\p\cuk(b)$. So, (PAL3) is fulfilled. Finally,
it is easy to verify (PAL4) and (PAL5). \sqs

\begin{pro}\label{catpal}
$\PAL$ is a category.
\end{pro}

\doc This follows immediately from \ref{pf1}(f), \ref{pf1}(h), \ref{compos} and \ref{compos1}.
\sqs

\begin{pro}\label{aliso}
Let $X$ be a locally compact Hausdorff space. Then the NCAs
$(RC(X),C_{\rho_X})$ and $(RC(\a X),\rho_{\a X})$ are
CA-isomorphic (see \ref{Alexprn} and \ref{stanlocn} for the
notations) and the maps $e_{X,\a X}$, $r_{X,\a X}$ are
CA-isomorphisms between them (see \ref{isombool} for the
notations).
\end{pro}

\doc By \ref{isombool}, we have only to show that $AC_{\rho_X} B$
iff $\cl_{\a X}(A)\rho_{\a X}\cl_{\a X}(B)$, for every $A,B\in
RC(X)$. This follows easily from the respective definitions.
Hence, the map $e_{X,\a X}:(RC(X),C_{\rho_X})\lra (RC(\a
X,\rho_{\a X})$ is a CA-isomorphism. Thus the map  $r_{X,\a X}$
 is also a CA-isomorphism. \sqs

\begin{theorem}\label{gendv}
The categories $\PLC$ and $\PAL$ are dually equivalent.
\end{theorem}

\doc We will define two contravariant functors
$$\Xi^a:\PAL\lra\PLC \mbox{ and }\Xi^t:\PLC\lra\PAL.$$\\
{\em I. The definition of}\ \  $\Xi^t$.

\medskip

For every $(X,\tau)\in\card{\PLC}$, we let
$\Xi^t(X,\tau)=\Psi^t(X,\tau)$ (see (\ref{psit1}) for $\Psi^t$).

Let $f:(X,\tau)\lra (Y,\tau\ap)\in \PLC(X,Y)$. We set
\begin{equation}\label{xit}
\Xi^t(f):\Xi^t(Y,\tau\ap)\lra\Xi^t(X,\tau), \ \ \
\Xi^t(f)(F)=cl_X(f\inv(\int_Y(F))).
\end{equation}
Put, for the sake of brevity, $\p_f=\Xi^t(f)$. We have to show
that $\p_f$ is a $\PAL$-morphism. Obviously, (PAL1) is fulfilled.
For verifying (PAL4), let $H\in CR(X)$. Then $f(H)$ is compact.
Since $Y$ is locally compact, there exists $F\in CR(Y)$ such that
$f(H)\sbe\int(F)$. Now we obtain that $H\sbe
f\inv(\int(F))\sbe\int(\cl(f\inv(\int(F))))=\int(\p_f(F))$, i.e.
$H\ll_{\rho_X}\p_f(F)$. Hence (PAL4) is checked.

Let now $F\in CR(Y)$. Then $\p_f(F)=\cl(f\inv(\int(F)))\sbe
f\inv(F)$. Since $f\inv(F)$ is compact (because $f$ is perfect),
$\p_f(F)\in CR(X)$. Therefore, (PAL5) is fulfilled.

By \ref{perfect1}, $f$ has a continuous extension $\a(f):\a X\lra
\a Y$. Set $\p_{\a f}=\Phi^t(\a(f))$ (see Theorem \ref{dvth} for
$\Phi^t$). Then, by Theorem \ref{dvth}, $\p_{\a f}$ is a
$\VAL$-morphism. We will prove that
\begin{equation}\label{comd}
r_{X,\a X}\circ\p_{\a f}=\p_f\circ r_{Y,\a Y}
\end{equation}
(see \ref{isombool} for the notations), i.e. that, for every $G\in
RC(Y)$, the following equality holds:
\begin{equation}\label{comdx}
X\cap\p_{\a f}(\cl_{\a Y}(G))=\p_f(G),
\end{equation}
 or, in other words, that
$$X\cap\cl_{\a X}((\a(f))\inv(\int_{\a Y}(\cl_{\a Y}(G))))=\cl_{
X}(f\inv(\int_Y(G))).$$
Since the last equality follows easily from the obvious inclusions
$\int_Y(G)\cup\{\iy\}\spe\int_{\a Y}(\cl_{\a Y}(G))\spe\int_Y(G)$,
(\ref{comd}) is proved. Therefore, $\p_f=r_{X,\a X}\circ\p_{\a
f}\circ e_{Y,\a Y}$ (see \ref{isombool}). Since $\p_{\a f}$
satisfies (DVAL2), we obtain that $\p_f$ satisfies (PAL2).

For establishing (PAL3), let $F\in CR(Y),G\in RC(Y)$ and
 $F\ll_{\rho_Y} G$. Then  $F\ll_{C_
 {\rho_Y}} G$ and hence, by
 \ref{aliso},  $F\ll_{\rho_{\a Y}} \cl_{\a Y}(G)$. Thus, (DVAL3)
 implies that
\begin{equation}\label{dve}
 (\p_{\a f}(F^{*\a}))^{*\a}\ll_{\rho_{\a X}}\p_{\a
 f}(\cl_{\a Y}(G)),
\end{equation}
  where $``{}^{*\a}$" is used as a common notation of  the
  complement
 in the Boole\-an algebras $RC(\a X)$ and
 $RC(\a Y)$. Since, for every $H\in
 RC(X)$, $X\cap(\cl_{\a X}(H))^{*\a}=r_{X,\a X}((\cl_{\a
 X}(H))^{*\a})=(r_{X,\a X}(\cl_{\a X}(H))^*=H^*$, we get, using
 again \ref{aliso}, that $(X\cap\p_{\a
 f}(F^{*\a}))^*\ll_{C_{\rho_X}} (X\cap\p_{\a f}(\cl_{\a Y}(G)))$;
 then,
 applying twice (\ref{comdx}),  the equality
 $F^{*\a}(=(e_{Y,\a Y}(F))^{*\a})=e_{Y,\a Y}(F^*)$ and (\ref{comd}), we obtain that
$(\p_f(F^*))^*\ll_{C_{\rho_X}} \p_f(G)$, i.e. (PAL3) is fulfilled.

Now, we will verify (PAL6). Let $F\in RC(Y)$; then $\cl_{\a
Y}(F)\in RC(\a Y)$ and hence, by (DVAL4),
$$\p_{\a f}(\cl_{\a Y}(F))=\bigvee\{\p_{\a f}(\cl_{\a Y}(G))\st
G\in RC(Y), \cl_{\a Y}(G)\ll_{\rho_{\a Y}}\cl_{\a Y}(F)\}.$$
 Since
$r_{X,\a X}$ is an isomorphism, we obtain that
$r_{X,\a X}(\p_{\a f}(\cl_{\a Y}(F)))= \bigvee\{r_{X,\a X}(\p_{\a
f}(\cl_{\a Y}(G)))\st G\in RC(Y), \cl_{\a Y}(G)\ll_{\rho_{\a
Y}}\cl_{\a Y}(F)\}.$
 Thus,
(\ref{comd}) and \ref{aliso} imply that
$\p_f(F)=\bigvee\{\p_f(G)\st G\in RC(Y), G\ll_{C_{\rho_Y}} F\}$.
So, (PAL6) is fulfilled.

 Therefore, $\p_f$ is
a $\PAL$-morphism.

Let $f\in\PLC(X,Y)$ and $g\in\PLC(Y,Z)$. We will prove that
$\Xi^t(g\circ f)=\Xi^t(f)\diamond\Xi^t(g)$. Put $h=g\circ f$,
$\p_h=\Xi^t(h)$, $\p_f=\Xi^t(f)$ and $\p_g=\Xi^t(g)$. Let
$\a(f):\a X\lra \a Y$, $\a(g):\a Y\lra\a Z$ and $\a(h):\a X\lra\a
Z$ be the continuous extensions of $f$, $g$ and $h$, respectively
(see \ref{perfect1}). Then, obviously, $\a(h)=\a(g)\circ\a(f)$.
Set $\p_{\a f}=\Phi^t(\a(f))$, $\p_{\a g}=\Phi^t(\a(g))$ and
$\p_{\a h}=\Phi^t(\a(h))$ Then, by Theorem \ref{dvth}, $\p_{\a
h}=(\p_{\a f}\circ\p_{\a g})\gek$. Now, using (\ref{comd}) and
\ref{isombool}, we get that $e_X\circ\p_h\circ r_Z=\p_{\a
h}=(e_X\circ\p_f\circ\p_g\circ r_Z)\gek$. Thus, for every $F\in
RC(\a Z)$, we have that
$\p_h(r_Z(F))=\bigvee\{(\p_f\circ\p_g)(r_Z(G))\st G\ll_{\rho_{\a
Z}} F\}$. Now, \ref{isombool} and \ref{aliso}  imply that
$\p_h=(\p_f\circ\p_g)\cuk$, i.e. $\p_h=\p_f\diamond\p_g$.

So, $\Xi^t:\PLC\lra\PAL$ is a contravariant functor.

\medskip

{\em II. The definition of}\ \ $\Xi^a$.

\medskip

For every $(A,\rho,\BBBB)\in\card{\PAL}$, we let
$\Xi^a(A,\rho,\BBBB)=\Psi^a(A,\rho,\BBBB)$ (see (\ref{phiapcn})
and (\ref{phiapc}) for $\Psi^a$).

Let $\p\in\PAL((A,\rho,\BBBB),(B,\eta,\BBBB\ap))$. We define the
map
$$\Xi^a(\p):\Xi^a(B,\eta,\BBBB\ap)\lra\Xi^a(A,\rho,\BBBB)$$
 by
the formula
\begin{equation}\label{xiap}
\Xi^a(\p)(\s\ap)=\{a\in A\st \mbox{ if } b\llcr a^*\mbox{ then
}(\p(b))^*\in\s\ap\},
\end{equation}
for every bounded cluster $\s\ap$ in $(B,C_\eta)$. Set, for the
sake of brevity, $\Xi^a(\p)=f_\p$,  $X=\Xi^a(A,\rho,\BBBB)$ and
$Y=\Xi^a(B,\eta,\BBBB\ap)$. We will show that $f_\p:Y\lra X$ is
well-defined and is a perfect map.

Let  $\p_C:(A,C_\rho)\lra(B,C_\eta)$ be defined by
$\p_C(a)=\p(a)$, for every $a\in A$. Then $\p_C$ is a
$\VAL$-morphism. Indeed, (DVAL3) follows from \ref{pf1}(c), and
the other three axioms are obviously fulfilled. Set
$f_\a=\Phi^a(\p_C)$. Then $f_\a:\a Y\lra\a X$ (see  Theorem
\ref{dvth} and (B1),  (B2) in the proof of Theorem \ref{roeperl}).
The definitions of $f_\p$ and $f_\a$ coincide on the bounded
clusters of $(B,C_\eta)$ (see (\ref{xiap}) and Theorem
\ref{dvth}); hence,  the right side of the formula (\ref{xiap})
defines a cluster in $(A,C_\rho)$ and $f_\a$ is an extension of
$f_\p$. Thus, if we show that $f_\a\inv(\ix)=\{\iy\}$, the map
$f_\p$ will be well-defined and will be a perfect map. Let us
prove that $f_\a(Y)\sbe X$, i.e. that if $\s\ap$ is a bounded
cluster in $(B,C_\eta)$ then $\s=f_\a(\s\ap)=f_\p(\s\ap)$ is a
bounded cluster in $(A,C_\rho)$. So, let $\s\ap$ be a bounded
cluster in $(B,C_\eta)$ and $\s=f_\a(\s\ap)$. Then \ref{bstar}
implies that there exists $b\in\BBBB\ap$ such that $b^*\nin\s\ap$.
By (PAL4), there exists $a\in\BBBB$ such that $b\le\p(a)$. Thus
$(\p(a))^*\le b^*$ and hence $(\p(a))^*\nin\s\ap$. By (BC1), there
exists $a_1\in\BBBB$ such that $a\llx a_1$. Then $a\llcr a_1$ and,
by the definition of $\s$, $a_1^*\nin\s$. Therefore
$a_1\in\BBBB\cap\s$, i.e. $\s$ is a bounded cluster in
$(A,C_\rho)$. Hence $f_\p(Y)=f_\a(Y)\sbe X$.
 Further, we have (by \ref{neogrn}) that
$\ix=A\stm \BBBB$ and $\iy=B\stm\BBBB\ap$. Let us show that
$f_\a(\iy)=\ix$. Set $\s\ap=\iy$ and $\s=f_\a(\s\ap)$. Let
$a\in\s$. Suppose that $a\in\BBBB$. Then, by (BC1), there exist
$a_1, a_2\in\BBBB$ such that $a\llx a_1\llx a_2$. Thus $a\llcr
a_1\llcr a_2$. Hence $a_1^*\llcr a^*$. Since $a\in\s$, the
definition of $\s$ implies that $(\p(a_1^*))^*\in\s\ap$. By
\ref{pf1}(c), we have that $(\p(a_1^*))^*\le\p(a_2)$. Therefore,
$\p(a_2)\in\s\ap$. Since $\p(a_2)\in\BBBB\ap$ (by (PAL5)), we
obtain a contradiction. Thus $\s\sbe A\stm\BBBB$. Now,
\ref{neogrn} and \ref{fact29} imply that $\s=A\stm\BBBB$, i.e.
$f_\a(\iy)=\ix$. Hence $f_\a\inv(X)=Y$. This shows that $f_\p$ is
a perfect map (because $f_\a$ is such). So, we have proved that
$f_\p\in\PLC(Y,X)$.

Let
$\p_i\in\PAL((A_i,\rho_i,\BBBB_i),(A_{i+1},\rho_{i+1},\BBBB_{i+1}))$
and $f_i=\Xi^a(\p_i)$ for $i=1,2$, $\p=\p_2\diamond\p_1$,
$f_\p=\Xi^a(\p)$ and $X_i=\Xi^a(A_i,\rho_i,\BBBB_i)$ for
$i=1,2,3$. We will prove that $f_\p=f_{1}\circ f_{2}$. Let
$(\p_i)_C:(A_i,C_{\rho_i})\lra (A_{i+1},C_{\rho_{i+1}})$ be
defined by $(\p_i)_C(a)=\p_i(a)$ for every $a\in A_i$, where
$i=1,2$. Then, as we know, $(\p_i)_C$ is a $\VAL$-morphism, for
$i=1,2$. Set $f_{i\a}=\Phi^a((\p_i)_C)$ for $i=1,2$,
$\psi=(\p_2)_C\,\ast\,(\p_1)_C$, $f_{\psi}=\Phi^a(\psi)$. Let
$\p_C:(A_1,C_{\rho_1})\lra (A_{3},C_{\rho_{3}})$ be defined by
$\p_C(a)=\p(a)$ for every $a\in A_1$. From the respective
definitions we obtain that, for every $a\in A_1$,
$\psi(a)=((\p_2)_C\circ(\p_1)_C)\gek(a)=(\p_2\circ\p_1)\cuk(a)=\p(a)$.
Thus, $\psi=\p_C$. Hence $f_{\psi}=\Phi^a(\p_C)$.
 We know
that $\Phi^a(A_i,C_{\rho_i})=\a X_i$, for $i=1,2,3$, and $f_{i\a}$
is a continuous extension of $f_i$, for $i=1,2$. The equality
$``\psi=\p_C$" implies that
 $f_{\psi}$ is a continuous extension of $f_\p$. From Theorem
 \ref{dvth} we get that $f_\psi=f_{1\a}\circ f_{2\a}$. Since
 $f_{1\a}\inv(X_1)=X_2$ and $f_{2\a}\inv(X_2)=X_3$, we conclude
 that $f_\p=f_{1}\circ f_{2}$.

We have proved that $\Xi^a:\PAL\lra\PLC$ is a contravariant
functor.

\medskip

{\em III. $\Xi^a\circ\Xi^t$ is naturally isomorphic to the
identity functor $Id_{\,\PLC}$.}

Recall that, for every $X\in\card{\PLC}$, the map
$t_X:X\lra(\Xi^a\circ\Xi^t)(X)$, where $t_X(x)=\s_x$ for every
$x\in X$, is a homeomorphism (see (\ref{homeo})). We will show
that $t^l:Id_{\,\PLC}\lra \Xi^a\circ\Xi^t$, where for every
$X\in\card\PLC$, $t^l(X)=t_X$, is a natural isomorphism.

Let $f\in\PLC(X,Y)$ and $f\ap=(\Xi^a\circ\Xi^t)(f)$,
$X\ap=(\Xi^a\circ\Xi^t)(X)$, $Y\ap=(\Xi^a\circ\Xi^t)(Y)$. We have
to prove that $t_Y\circ f=f\ap\circ t_X$. Let $\a(f):\a X\lra\a Y$
and $\a(f\ap):\a X\ap\lra\a Y\ap$ be the continuous extensions of
$f$ and $f\ap$, respectively (see \ref{perfect1}). Then, by
Theorem \ref{dvth}, we have that $t_{\a Y}\circ\a(f)=\a(f\ap)\circ
t_{\a X}$. Obviously, $t_{\a X}(\ix)=\{\cl_{\a X}(F)\st F\in
RC(X), \ix\in\cl_{\a X}(F)\}=\{\cl_{\a X}(F)\st F\in RC(X)\stm
CR(X)\}=\s_\infty^{(RC(\a X),\rho_{\a X})}=\infty_{X\ap}$, and,
analogously, $t_{\a Y}(\iy)=\infty_{Y\ap}$. Using \ref{isombool}
and taking the restrictions on $X$, we obtain that
 $t_Y\circ f=f\ap\circ t_X$, i.e.
 $Id_{\,\PLC}\cong\Xi^a\circ\Xi^t$.

\medskip

{\em IV. $\Xi^t\circ\Xi^a$ is naturally isomorphic to the identity
functor $Id_{\, \PAL}$.}

Recall that for every $(A,\rho,\BBBB)\in\card{\PAL}$, the function
$$\l_A^g:(A,\rho,\BBBB)\lra(\Xi^t\circ\Xi^a)(A,\rho,\BBBB)$$
 is an
LCA-isomorphism (see (\ref{hapisomn})). We will show that
$\l^g:Id_{\,\PAL}\lra \Xi^t\circ\Xi^a$, where for every
$(A,\rho,\BBBB)\in\card\PAL$, $\l^g(A,\rho,\BBBB)=\l_A^g$, is a
natural isomorphism.

Let $\p\in\PAL((A,\rho,\BBBB),(B,\eta,\BBBB\ap))$ and
$\p\ap=(\Xi^t\circ\Xi^a)(\p)$, $X=\Xi^a(A,\rho,\BBBB)$,
$Y=\Xi^a(B,\eta,\BBBB\ap)$. We have to prove that
$\l_B^g\diamond\p=\p\ap\diamond\l_A^g$. According to (\ref{diam})
and (\ref{cukf}), it is enough to show that
$\l_B^g\circ\p=\p\ap\circ\l_A^g$. Set $f=\Xi^a(\p)$. Hence
$\p\ap=\Xi^t(f)$. Let $\p_C:(A,C_{\rho})\lra (B,C_{\eta})$ be
defined by $\p_C(a)=\p(a)$ for every $a\in A$, and let $(\p\ap)_C$
be defined analogously. Then $\p_C$ and $(\p\ap)_C$  are
$\VAL$-morphisms. Set $f_\a=\Phi^a(\p_C)$ and
$(\p_C)\ap=\Phi^t(f_\a)$. We know that $f_\a:\a Y\lra \a X$ is a
continuous extension of $f$. By the proof of Theorem \ref{dvth},
$\l_B\circ\p_C=(\p_C)\ap\circ\l_A$ (see (\ref{h}) for $\l_A$ and
$\l_B$). Note that $\l_A:(A,C_\rho)\lra(RC(\a X),\rho_{\a X})$ and
$\l_B:(B,C_\eta)\lra(RC(\a Y),\rho_{\a Y})$. Let
$(\p\ap)_C:(RC(X), C_{\rho_X})\lra (RC(Y), C_{\rho_Y})$ be defined
by $(\p\ap)_C(F)=\p\ap(F)$, for every $F\in RC(X)$. Then, by
(\ref{comd}), $(\p\ap)_C\circ r_X=r_Y\circ(\p_C)\ap$. By
(\ref{lbg}), $r_X\circ\l_A=\l_A^g$ and $r_Y\circ\l_B=\l_B^g$. The
last three equalities imply that $\l_B^g\circ\p=\p\ap\circ\l_A^g$.
 Thus $Id_{\,\PAL}\cong\Xi^t\circ\Xi^a$.
\sqs

\begin{theorem}\label{phiinj}
Let $\p$ be a $\PAL$-morphism. Then $\p$ is an injection iff\/
$\Xi^a(\p)$ is a surjection.
\end{theorem}

\doc Let $\p\in\PAL((A,\rho,\BBBB),(B,\eta,\BBBB))$ and let
$\p_C:(A,C_\rho)\lra(B,C_\eta)$ be defined by the formula
$\p_C(a)=\p(a)$, for every $a\in A$. Then $\p_C$ is a
$\VAL$-morphism. Setting  $f=\Xi^a(\p)$, we obtain that
$\a(f)=\Phi^a(\p_C)$ (see the proof of Theorem \ref{gendv}).
Obviously, $\a(f)$ is a surjection iff $f$ is a surjection. By a
theorem of de Vries (\cite[Theorem 1.7.1]{dV}), $\Phi^a(\p_C)$ is
a surjection iff $\p_C$ is an injection. Hence, $f$ is a
surjection iff $\p$ is an injection. \sqs

It is clear that if we want to build $\PAL$ as a category dually
equivalent to the category $\PLC$ then the axiom (PAL5) is
indispensable for describing the morphisms of the category $\PAL$.
With the next simple example we show that the axiom (PAL4) cannot
be dropped as well.

\begin{exa}\label{pal4e}
\rm Let $(A,\rho,\BBBB)$ be a CLCA and $\BBBB\neq A$. Then
$(A,\rho_s,A)$ is also a CLCA (by \ref{extrcr}). Obviously, the
map $i:(A,\rho,\BBBB)\lra (A,\rho_s,A)$, where $i(a)=a$, for every
$a\in A$,  satisfies the axioms (PAL1)-(PAL3), (PAL5), (PAL6) but
it does not satisfy the axiom (PAL4). If we suppose that our
duality theorem is true without the presence of the axiom (PAL4)
in the definition of the category $\PAL$ then we will obtain, by
Theorem \ref{phiinj}, that there exists a continuous map from a
compact Hausdorff space onto a locally compact non-compact
Hausdorff space, a contradiction.
\end{exa}

\begin{fact}\label{lcase}
For every LCA $(A,\rho,\BBBB)$, the triple $(A,\rho_s,\BBBB)$ is
also an LCA (see \ref{extrcr} for $\rho_s$); if $(A,\rho,\BBBB)$
is a CLCA then the map $i:(A,\rho,\BBBB)\lra (A,\rho_s,\BBBB)$,
where $i(a)=a$, for every $a\in A$, is a $\PAL$-morphism.
\end{fact}

\doc Since $a\ll_{\rho_s} a$, for every $a\in A$, the axiom (BC1)
of \ref{locono} is clearly fulfilled. Obviously, for every $a,b\in
A$, $a\ll_\rho b$ implies $a\ll_{\rho_s} b$. This implies that the
axiom (BC3) is also satisfied. For checking (BC2), let $a,b\in A$
and $a\rho_s b$. Then $a\wedge b\not= 0$. Since $b=\bigvee\{c\st
c\in \BBBB, c\ll_\rho b\}$, we have that $b=\bigvee\{c\st c\in
\BBBB, c\wedge b^*=0\}$. Hence $a\wedge b=\bigvee\{a\wedge c\st
c\in \BBBB, c\wedge b^*=0\}$. Thus, there exists $c\in\BBBB$ such
that $c\wedge b^*=0$ and $a\wedge c\not= 0$. Therefore, there
exists $c\in\BBBB$ such that $a\rho_s(c\wedge b)$. So,
$(A,\rho_s,\BBBB)$ is  an LCA. The rest is clear.  \sqs

Recall that a topological space $X$ is said to be {\em extremally
disconnected}\/ if for every open set $U\sbe X$, the closure
$\cl_X(U)$ is open in $X$. Clearly, a topological space $X$ is
 extremally disconnected iff $RC(X)$ consists only of clopen sets.

\begin{pro}\label{comextr}
Let $(A,\rho,\BBBB)$ be a CLCA and $X=\Xi^a(A,\rho,\BBBB)$.
Then:\\
(a)(\cite{R}) $X$ is a compact Hausdorff space iff\/ $\BBBB=A$;\\
(b) $X$ is an extremally disconnected locally compact Hausdorff
space iff $\rho=\rho_s$ (see \ref{extrcr} for $\rho_s$).
\end{pro}

\doc The assertion (a) is obvious.

\noindent(b) Recall that, by (\ref{hapisomn}),
$\l_A^g:(A,\rho,\BBBB)\lra (RC(X),\rho_X)$ is an LCA-isomorphism.

Let $X$ be extremally disconnected. Then, for every $a,b\in A$,
$\l_A^g(a\wedge
b)=\l_A^g(a)\wedge\l_A^g(b)=\cl(\int(\l_A^g(a)\cap\l_A^g(b)))=
\l_A^g(a)\cap\l_A^g(b)$. Hence $a\wedge b\not= 0$ iff $a\rho b$.
Thus $\rho=\rho_s$ (see \ref{extrcr}).

Conversely, let $\rho=\rho_s$. Then, for every $a\in A$,
$a\ll_{\rho} a$. Since for every $a,b\in A$,  $a\ll_\rho b$ iff
$\l_A^g(a)\sbe\int_X(\l_A^g(b))$, we get that for every $a\in A$,
$\l_A^g(a)\sbe\int_X(\l_A^g(a))$, i.e. $\l_A^g(a)$ is a clopen
set. Therefore, $X$ is extremally disconnected. \sqs

Note that from \ref{comextr}(b), \ref{lcase} and Theorem
\ref{phiinj}, we obtain immediately an easy proof of the following
well-known fact: every locally compact Hausdorff space $X$ is a
perfect image of an extremally disconnected locally compact
Hausdorff space $Y$.

\begin{theorem}\label{clemb}
Let $X$ and $Y$ be two locally compact Hausdorff spaces,
$\Xi^t(X)=(A,\rho,\BBBB)$ and\/ $\Xi^t(Y)=(B,\eta,\BBBB\ap)$. Then
a map $f:X\lra Y$  is a closed embedding iff the map $\p=\Xi^t(f)$
satisfies the following two conditions:\\
\noindent(1) $\fa a,b\in A$ with $a\llcr b$ there exists
$c\in B$ such that $a\llcr \p(c)\llcr b$;\\
\noindent(2) $\fa a,b\in B$, $\p(a)\llcr\p(b)$ iff there exist
$a_1, b_1\in B$ such that $a_1\llce b_1$ and $\p(a_1)=\p(a)$,
$\p(b_1)=\p(b)$.
\end{theorem}

\doc Obviously, $f:X\lra Y$ is a closed embedding iff the map $\a(f):\a
X\lra\a Y$ is an embedding (note that every closed embedding is a
perfect map and see \ref{perfect1} for $\a(f)$).  De  Vries proved
(see \cite[Theorem 1.7.3]{dV}) that $\a(f)$ is an embedding iff
the following two conditions are sa\-tisfied: (a) for every
$F,G\in RC(\a X)$ with $F\ll_{\rho_{\a X}} G$, there exists
$H\in\Phi^t(\a(f))(RC(\a Y))$ such that $F\ll_{\rho_{\a X}}
H\ll_{\rho_{\a X}} G$, and (b) for every $F,G\in RC(\a Y)$,
$\Phi^t(\a(f))(F)\ll_{\rho_{\a X}} \Phi^t(\a(f))(G)$ iff there
exist $F_1,G_1\in RC(\a Y)$ such that $F_1\ll_{\rho_{\a Y}} G_1$
and $\Phi^t(\a(f))(F_1)=\Phi^t(\a(f))(F)$,
$\Phi^t(\a(f))(G_1)=\Phi^t(\a(f))(G)$. Now, using \ref{aliso} and
(\ref{comd}), it is easy to obtain that $f$ is a closed embedding
iff $\p$ satisfies conditions (1) and (2). \sqs

\begin{notas}\label{notac}
\rm Let us denote by $\PLCC$ the full subcategory of the category
$\PLC$ whose objects are all connected locally compact Hausdorff
spaces. Let $\PALC$ be the full subcategory of the category $\PAL$
whose objects are all connected CLCAs.
\end{notas}

\begin{theorem}\label{gencdv}
The categories $\PLCC$ and $\PALC$ are dually equivalent.
\end{theorem}

\doc It follows immediately from Theorem \ref{gendv} and Fact \ref{confact}.


\baselineskip = 0.75\normalbaselineskip

\baselineskip = 1.00\normalbaselineskip

\vspace{0.25cm}

 Faculty  of Mathematics and Informatics,

Sofia University,

5 J. Bourchier Blvd.,

1164 Sofia,

Bulgaria

\end{document}